\documentstyle[12pt]{article} 
\newcommand{\Uq}{U_q(osp(1|2))} 
 
\newcommand{\lam}{\lambda} 
 
\newcommand{\x}{\otimes}

\newcommand{\ba}{\begin{eqnarray}} 
\newcommand{\na}{\end{eqnarray}} 
\newcommand{\ban}{\begin{eqnarray*}} 
\newcommand{\nan}{\end{eqnarray*}}

\newcommand{\Zz}{{\bf Z}_2} 
\newcommand{\Zn}{{\bf Z}_N}

\begin{document} 
\title{\small\bf  LICKORISH INVARIANT AND QUANTUM OSP(1$\mid$2)}  
\author{ \small R. B. ZHANG$^\dagger$ \ \ and \ \ H. C. LEE$^*$ \\ 
\small $^\dagger$Department of Pure Mathematics\\ 
\small University of Adelaide, Adelaide, S. A. 5001, Australia\\ 
\small $^*$Physics Department\\ 
\small National Central University, Chungli, Taiwan, R. O. China}  
\date{ }
\maketitle  

\vspace{3cm}
\noindent 
Lickorish's method for constructing topological invariants of 
$3$ - manifolds is generalized to the quantum supergroup setting.   
An invariant is obtained by applying this method to  
the Kauffman polynomial arising from  the vector representation of 
$U_q(osp(1|2))$. A transparent proof is also given showing 
that this invariant is equivalent to the $U_q(osp(1|2))$ invariant 
obtained in an earlier publication.

\vspace{3cm} 
\noindent {\bf Introduction}

Since Jones' seminal work\cite{Jones}, the theories of knots and 
$3$ - manifolds have made dramatical progress ( See \cite{Turaev} 
for a review).  By now several approachs are available for 
constructing the so called quantum invariants of $3$ - manifolds, 
notably, the quantum field theoretical approach\cite{Witten}, 
the quantum group approach\cite{Reshetikhin-Turaev}, 
Lickorish's recoupling theory\cite{cable}, 
the $6 j$ symbol method of Turaev - Viro\cite{Viro},  
and the conformal field theoretical method\cite{Kohno}.   
All these approachs originated from theoretical physics, 
but they differ from one another very significantly 
in the mathematical formulations, and  each method has its own advantages 
in addressing specific problems. Therefore, it is important to 
further develop the different approaches, even though it is believed  
that they all give rise to the same invariants of $3$ - manifolds
(The Turaev - Viro invariant is known to be the square of the 
norm of the Reshetikhin-Turaev invariant). 

The recoupling theory was first introduced in \cite{cable} by 
Lickorish, who used the representation theory of the Temperley - Lieb 
algebra to reproduce the Jones - Witten - Reshetikhin - Turaev 
invariants.  Since then the method has been developed extensively
by other people\cite{Wenzl}\cite{Vogel}\cite{Kauffman}. 
In this letter, we aim to 
extend the recoupling theory in another direction, namely, to 
incorporate supersymmetry. We will also investigate the connection 
of the recoupling theory with the Reshetikhin - Turaev formalism 
in the quantum supergroup setting.  
    
Parallel to the quantum group formulation of topological invariants 
of links and $3$ - manifolds, there exists a supersymmetric 
version\cite{Zhang} 
based on the theory of Lie superalgebras and quantum supergroups. 
Due to the vast difference between the representation theory of 
quantum supergroups and that of the ordinary quantum groups, 
the associated topological invariants in both cases also exhibit 
different features. For example, there exist infinitely many families 
of multi - parameter generalizations of the Alexander - Conway invariant,
arising from the so - called typical irreps.  Such invariants can not 
be obtained within the framework of ordinary quantum groups. 
In this letter, we will apply the recoupling theory to the Kauffman 
polynomial associated with the vector representation of the quantum 
supergroup $\Uq$  to construct the corresponding topological invariant of 
$3$ - manifolds.  We will also provide a transparent proof showing 
that the resultant invariant is equivalent to that constructed 
in \cite{osp}.  In doing so, we establish  a precise relationship 
between the recoupling theory and Reshetikhin - Turaev method within 
the context of our study. 

Although we have limited ourselves to the quantum supergroup $\Uq$ 
here, the method developed for constructing the Lickorish invariant 
can be readily applied to self - dual atypical irreps of any quantum 
supergroup. It also appears to be possible to extend the formalism 
to include typical irreps, which of course are much more interesting 
due to their connections with the generalized Alexander - Conway 
polynomials, and possibly generalized versions of the Casson 
invariant\cite{Akbulut}. 
Results on this problem will be reported in a separate publication.

\bigskip 
\noindent
{\bf $\Uq$}

We will work on the complex field $\bf C$. 
Let $q$ be an $N$ - th primitive root of unity with $N$ a
positive odd integer satisfying    
\ban 
N=2r+1, &0<r\in{\bf Z}_+. 
\nan 
The quantum supergroup $\Uq$ is a $\Zz$ graded Hopf algebra,  
with the underlying algebra generated 
by $\{e, f, K^{\pm 1}\}$ subject to the relations 
\ba 
e f + f e &=& {{K- K^{-1}}\over{q-q^{-1}}}, \nonumber \\ 
K e K^{-1}&= &q e, \nonumber \\
K f K^{-1}&=& q^{-1}f,\nonumber \\ 
e^{2N}=f^{2N}&=&0,\nonumber \\
 K^{\pm N}&=&1,  
\na 
where the elements $e$ and $f$ are odd, while $K^{\pm 1}$ are even. 
The co - multiplication $\Delta: \Uq\rightarrow\Uq\x\Uq$ is 
given by 
\ban 
\Delta(e)&=&e\x K + 1\x e,\\ 
\Delta(f)&=&f\x 1 + K^{-1}\x f,\\ 
\Delta(K^{\pm})&=&
K^{\pm1}\x K^{\pm1}.  
\nan  

It is well known that the $\Uq$ so defined has the structures of 
a quasi triangular $\Zz$ graded Hopf algebra. We denote by $R$ 
its universal $R$ matrix and set $R=\sum \alpha_t\otimes \beta_t$. 
Then $R$ provides an isomorphism between the two algebras $\Delta[\Uq]$ 
and $\Delta'[\Uq]$, where $\Delta'$ is the opposite co - multiplication. 
Furthermore, $R$  also satisfies the quantum Yang - Baxter equation. 
Define $ v= K^{-1}\sum_t (-1)^{[\alpha_t]}S(\beta_t)\alpha_t$, 
where $S$ is the antipode of $\Uq$. Then $v$ belongs to the center 
of $\Uq$. 

Any $\Uq$ module $W$ is a $\Zz$ graded vector space $W=W_0\oplus W_1$, 
where $W_0$ and $W_1$ are the even and odd subspaces respectively. 
Let $W'$ be the $\Zz$ graded vector space with $W'_0=W_1$, 
and $W'_1=W_0$. Then $W'$ has a natural $\Uq$ module structure. 
These two modules are evidently isomorphic, with the isomorphism given
by a homogeneous degree $1$ linear mapping.
Let $f: W\rightarrow W$ be any module homomorphism. Then $f$ gives 
rise to a module homomorphism $f: W'\rightarrow W'$ in a natural way. 
Now  the $q$ - superdimensions of $f$ satisfy the following   obvious 
relation 
\ban 
Str_W (K f) &=& - Str_{W'} (K f). 
\nan

In \cite{osp}, we classified all the irreducible representations 
of this quantum supergroup. It was shown that there exist only a finite 
number of irreducible representations, which are all finite dimensional, 
and each of them is  uniquely characterized by an integer in 
$\Zn =\{ 0, 1, ..., N-1\}$.   
Let $V^+(\lam)$ be an irreducible $\Uq$ module, 
then it possesses a unique highest weight vector $v_0(\lam)$, which 
is assumed to be even,  such that 
$e v_0(\lam)=0$, $K v_0(\lam)=q^\lam v_0(\lam)$, $\lam\in \Zn$. 
A basis for $V^+(\lam)$ is given by 
$\{ v_0(\lam), v_1(\lam), ..., v_{2\lam}(\lam)\}$, 
where
$v_{i+1}(\lam)=f v_i(\lam)$, $f v_{2\lam}(\lam)=0$. 
We will denote by $V^-(\lam)$ the $\Uq$ module isomorphic to $V^+(\lam)$
but with an odd highest weight vector.

The $q$ - superdimension of $V^+(\lam)$ is given by 
\ban SD_q(\lam)&:=& Str_{V^+(\lam)}\left( K\right)\\ 
          &=& { {q^{\lam+1}+ q^{-\lam}}\over{q+ 1}}. \nan     
An important fact is that all the irreps have nonzero $q$ - 
superdimensions, and for this reason, all the irreps had to be 
included in constructing  the $3$ - manifold invariant of \cite{Zhang}. 
However, the $S$ - matrix arising from the Hopf link is singular, 
thus $\Uq$ does not qualify as a $\Zz$ graded modular Hopf algebra.

The smallest nontrivial irreducible $\Uq$ module  $V^+(1)$ 
will play an important role in the remainder of the note. 
We will denote it by $V$, and the associated irrep by $\pi$. 
It is assumed through out the paper that the highest weight vector
of $V$ is {\em even}.  It is important to  
observe that the irrep is self dual, that is, there exists a 
homogeneous degree zero isomorphism between $V$ and the dual module 
$V^*$.  The tensor product module $V\otimes V$ decomposes into 
\ban 
V\otimes V&=& V^+(2)\oplus V^-(1) \oplus V^+(0).   
\nan 
Thus the braid generator $\check{R}= P( \pi\otimes\pi )R$ 
satisfies the following third order polynomial relation 
\ban 
(\check{R} - q) ( \check{R} + q^{-1}) (\check{R} - q^{-2})&=&0. 
\nan 
This in particular implies that the link invariant arises from 
$V$ is the Kauffman polynomial.

Consider the central element of $\Uq$ defined by  
\ban 
{\hat C}^{(k)}&=&Str_{V^{\x k}}[(\pi^{\x k} \Delta^{(k-1)}\x id)
( v^{-1} K^{-1}\x 1)R^T R], \ \ \  k \in {\bf Z}_+. 
\nan 
Acting on the irreducible $\Uq$  module $V^+ (\lam)$, 
${\hat C}^{(k)}$ takes the eigenvalue 
\ban 
\chi_\lam ({\hat C}^{(k)})&=& Str_{V^{\x k}}[\pi^{\x k} \Delta^{(k-1)}
( v^{-1} K^{-2\lam -1})]. 
\nan 
Note that $\chi_\lam ({\hat C}^{(k)})$ is a finite sum of powers of $q$, 
thus it is consistent to first evaluate $\chi_\lam ({\hat C}^{(k)})$ at 
{\em generic} $q$ then specialize it to the $N$ - th root of unity. 
This way we obtain 
\ban 
\chi_\lam ({\hat C}^{(k)})&=&
\sum_{j=0}^k b_j ^{(k)} { {q^{(j+1)(2\lam +1)} + q^{-j(2\lam +1)}} 
\over{ 1 + q^{2\lam +1}} } q^{j (j+1)}, 
\nan 
where $b_j^{(k)}$ are a set of complex numbers determined by the
following recursion relations 
\ba
b_j^{(n+1)}&=& b_{j-1}^{(n)} - b_{j}^{(n)}  + b_{j+1}^{(n)}, 
           \  \  j>0\nonumber \\ 
b_0^{(n+1)}&=& b_1^{(n)},    \label{B} 
\na 
with the initial condition $b_0^{(0)}=1$, $b_j^{(0)}=0$, 
$\forall j>0$. It is easy to see that 
\ban 
b^{(n)}_n=1,   & b^{(n)}_{n+j}=0, \ \ \ \forall j>0.\\  
\nan

\noindent 
{\bf Lickorish Invariant} 

We construct the Lickorish invariant of $3$ - manifolds in this section. 
Although we only consider the quantum supergroup $\Uq$, the method 
developed here is general, and can be readily applied to self - dual 
atypical irreps of any other quantum supergroups.

We will need some facts about the quantum supergroup approach to 
invariants of framed links, which we recall here. 
The Reshetikhin - Turaev approach to link invariants was generalized 
to quantum supergroups in \cite{Zhang},  where a functor from the category 
of coloured ribbon graphs to the category of representations of 
$\Zz$ graded ribbon Hopf algebras was constructed. 
In plain term, this functor associates each coloured ribbon graph 
with a homomorphism of $\Zz$ graded modules of a quantum supergroup, 
where the modules are associated with the `colour' of the graph. 
In particular, an over crossing is represented by the universal $R$ 
matrix, and  an under crossing by $R^{-1}$. The composition of 
coloured ribbon graphs corresponds to the composition of homomorphisms 
of quantum supergroup modules, and the juxtaposition of ribbon 
graphs to tensor product of module homomorphisms.  The precise 
definition was given in \cite{Zhang} in explicit form, 
and we refer to that paper for details. 
Here we merely discuss a few aspects of the functor, 
which will be used extensively later.

Given any oriented $(k, l)$ ribbon graph $\Gamma$, 
we colour each of its components by the vector module $V$ of $\Uq$. 
The Reshetikhin - Turaev functor maps the coloured ribbon graph to 
a module homomorphism $F(\Gamma_V): 
V^{\x k}\rightarrow V^{\x l}$. Now reversing the orientation of any 
component of the ribbon graph will change the colouring module for that 
component to its dual module.   As we have already  pointed out, 
$V$ is self - dual, thus $F(\Gamma_V)$ is independent of the orientation of 
$\Gamma$.  Also,  on any $(0, 0)$ ribbon graph, this module homomorphism 
yields the Kauffman polynomial for ribbon graphs. 

Consider the $(k, k)$ ribbon graphs given in Figure 1 and Figure 2 
respectively,  

\vspace{6cm}  
\begin{center} 
Figure 1 \hspace{5cm}  Figure 2  
\end{center} 

\noindent 
where Fig. 1 has $n$ annuli. We colour the ribbons of Fig. 2 
from left to right by the modules $W_1, ..., W_k$ respectively. 
We also colour the ribbons of Fig. 1 in the same way, 
but colour each of its annuli by $V$. 
Set $W= W_1\x ...\x W_k$, and denote the resultant 
coloured coloured ribbon graphs of Fig. 1 and Fig. 2 by 
$\phi^{(n)}_W$ and $\zeta_W$ respectively. 
Then the functor $F$ gives 
\ban
F(\phi^{(n)}_W)&=&{\hat C}^{(n)}: W\rightarrow W,\nonumber \\
F(\zeta_W)&=&v: W\rightarrow W,
\nan
where the central elements ${\hat C}^{(n)}$ and $v$ act on 
the tensor product module $W$ via the co - multiplication 
$\Delta^{(k-1)}$.

Let us now construct the Lickorish invariant. 
Lickorish's construction uses two fundamental theorems from the 
theory of $3$ - manifolds. In the earlier 1960s,  
Lickorish and  Wallace proved 
that each framed link in $S^3$ determines 
a compact, closed, oriented $3$ - manifold, 
and every such $3$ - manifold is obtainable by surgery 
along a framed link in $S^3$\cite{Lickorish}.
Further advances along this line were obtained by Kirby, Craggs, 
and Fenn and Rourke\cite{Kirby}\cite{Fenn}, 
who  proved that
orientation preserving homeomorphism classes of
compact, closed, oriented $3$ - manifolds correspond bijectively 
to equivalence classes of framed links in $S^3$, where 
the equivalence relation is generated by the Kirby moves.

Let $L$ be a framed link in $S^3$ with $m$ components $L_1, ..., 
L_m$.  We arbitrarily assign an orientation to each of its components. 
The resultant oriented framed link can be represented in a unique way 
by an oriented ribbon graph $\Gamma(L)$ \cite{Zhang}. Now we consider the 
oriented ribbon graph $\Gamma(L^{\{l_1, ..., l_m\}})$ obtained
from $\Gamma(L)$ in the following way: 
we replace each $L_i$ of $L$ by a cable of $l_i$ copies of $L_i$ with 
the same orientation, where $l_i\in\{0, 1, ..., N-1\}$. This leads to 
an oriented framed link with $\sum_{1\le i\le m} l_i$ components,  
the associated oriented ribbon graph of which is 
$\Gamma(L^{\{l_1, ..., l_m\}})$.              
We colour each component of $\Gamma(L^{\{l_1, ..., l_m\}})$ by the 
vector module $V$ of $\Uq$, and denote the corresponding coloured 
ribbon graph by $\Gamma(L^{\{l_1, ..., l_m\}})_V$. Applying the 
Reshetikhin - Turaev functor $F$ to it leads to a complex number 
$F(\Gamma(L^{\{l_1, ..., l_m\}})_V)$. 
The self duality of $V$ implies that this  number 
is independent of the arbitrarily chosen orientation of $L$. 

Construct 
\ba 
\Sigma (L) &=& \sum_{l_1, ..., l_m=0}^{N-1}\  \prod_{i=1}^m d^{(l_i)}
\ F(\Gamma(L^{\{l_1, ..., l_m\}})_V), \label{Sigma} 
\na 
where the $d^{(l_i)}$ are a set of constants chosen in such a way 
that $\Sigma (L)$ is invariant under the positive Kirby moves. 
Needless to say, the critical problem is whether such $d^{(l_i)}$
exist.  We will study this problem at great length later. Here 
let us take as granted their existence under 
the further assumption that 
\ba 
z&:=& \sum_{l=0}^{N-1} d^{(l)} F(\Gamma({\cal O}_{-1}^{\{l\}})_V) \label{z} \\ 
&\ne & 0, \nonumber 
\na 
where ${\cal O}_{-1}$ represents the unknot with framing number 
$-1$, and ${\cal O}_{-1}^{\{l\}}$ the framed link obtained by 
extending the framed unknot to $l$ parallel copies. Then 
\ba 
{\cal F}(M_L)&=& z^{-\sigma(A_L)} \Sigma (L), 
\na 
is a topological invariant of the $3$ - manifold $M_L$ obtained 
by surgery along the framed link $L$. Here $A_L$ is the linking matrix 
of $L$, and $\sigma(A_L)$ is the number of nonpositive eigenvalues 
of $A_L$.     
 
Assuming the properties of $d^{(l)}$, we can easily show that 
${\cal F}(M_L)$ is indeed a topological invariant of $M_L$:  
under the positive Kirby moves, $\sigma(A_L)$ is not changed,  
thus ${\cal F}(M_L)$ is invariant. Let $L'$ denote the split link  
$L\cup {\cal O}_{-1}$. Then $\Sigma (L') = z \Sigma (L)$. 
Since $\sigma(A_L') = \sigma(A_L) + 1$, 
${\cal F}(M_L)$ is invariant under the special negative Kirby 
move as well, and the proof is completed. 

Let us now construct the $d$'s.  In Lickorish's original paper, 
the Jones polynomial was used for constructing the Witten - 
Reshetikhin - Turaev invariant. There the Jones - Wenzl theory 
of the Temperley - Lieb algebra played an important role in 
the determination of the constants $d^{(l)}$. In our case, the 
Lickorish invariant of $3$ - manifolds is built from the Kauffman 
polynomial, which has a deep connection with the Birman - Wenzl 
algebra. It should be possible to obtain the $d^{(l)}$ by using 
only the representation theory of this algebra. However, 
such a method will not provide us with much information on the 
relationship between the Lickorish approach and the Reshetikhin - Turaev 
construction of $3$ - manifold invariants. 
On the other hand, the representation theory of $\Uq$ affords a 
common basis for both approaches, and also provides a  
much more powerful tool for determining the $d$'s. 

Consider the module homomorphisms associated with the coloured 
ribbon graphs of Figure 1 and Figure 2.  Let $f:  W \rightarrow W$ 
be any $\Uq$ module homomorphism.  Then the vanishing of the following 
$q$ - supertrace 
\ba 
Str_W\left[ K f \left(\sum_{l=0}^{N-1} d^{(l)} F(\phi^{(l)}_W) 
- F(\zeta_W)\right)\right]&=&0,   \label{master} 
\na 
for all $k$, arbitrary $W_i$'s and $f$, will guarantee the invariance 
of $\Sigma (L)$ under the positive Kirby moves.  A sufficient condition 
for equation (\ref{master}) to hold is that the central 
element of $\Uq$ defined by    
\ba 
\delta &=& \sum_{l=0}^{N-1} d^{(l)} {\hat C}^{(l)} - v,  
\na 
takes $0$ eigenvalue in all irreducible representations of $\Uq$, 
that is, 
\ba 
\sum_{l=0}^{N-1} d^{(l)} \chi_\lam ({\hat C}^{(l)})
&=& q^{-\lam(\lam+1)}, \ \ \ \ \forall \lam\in{\bf Z}_N. \label{equation} 
\na

Therefore, the problem of constructing $\Sigma(L)$ is now reduced to 
that of solving equation (\ref{equation}).  To do this, we 
introduce the $N\times N$ matrix $B = (b_{\mu \nu})_{\mu, \nu=0}^{N-1}$ 
with the entries given by $b_{\mu \nu} = b^{(\mu)}_\nu$, where 
$b^{(\mu)}_\nu$ are defined by the relations (\ref{B}). Write
$d=(d^{(0)}, d^{(1)},..., d^{(N-1)})$, and define      
\ban  
b &=&d B. 
\nan 
Since $B$ is lower triangular with all diagonal elements being $1$, 
there exists a unique $d$ corresponding to any given $b$.  Now in 
terms of the components of $b$,  equation (\ref{equation}) can be 
rewritten as 
\ba 
\sum_{\mu=0}^{N-1} b_\mu {{q^{(\mu+1)(2\lam +1)} + q^{-\mu(2\lam +1)}}
\over{ 1 + q^{2\lam +1}} } q^{\mu(\mu+1)}&=& q^{-\lam(\lam+1)}, 
\ \ \ \ \lam=0, 1, ..., N-1.  
\na  
It is precisely this equation which appeared in the Reshetikhin - 
Turaev construction of $3$ - manifold invariants using the quantum 
supergroup $\Uq$. The most general solution of this equation was obtained 
in \cite{osp}, which we quote  below
\ba
b_\mu&=&{{(1+q) q^{\left({{N+1}\over{2}}\right)^2} 
G_{-1} SD_q(\mu)}\over{2 N}} 
+x_\mu - x_{N-\mu -1},   \ \ \ \mu=0, 1, ..., N-1, 
\na 
where $x_\mu$ are arbitrary complex parameters,    
$G_{-1}=\sum_{\lam\in\Zn} q^{-\lam^2}$, and $SD_q(\mu)$ 
denotes the $q$ - superdimension of $V^+(\mu)$.

To compute $z$, we apply the formula  
\ban 
F(\Gamma({\cal O}_{-1}^{\{l\}})_V)&=& Str_{V^{\x l}} \left( v K\right) \\ 
&=&\sum_{j=0}^l b_j^{(l)}\  { { q^{j+1} + q^{-j}}\over{1+q}}\ q^{-j (j+1)},  
\nan 
to cast it into the form  
\ba
z&=&\sum_{\lam=0}^{N-1}b_\lam \ q^{-\lam (\lam +1)}\  SD_q(\lam),
\na
which coincides with the quantity also denoted by $z$ in \cite{osp}.  
Now $ z=q^{{N+3}\over{2}} \left( { {G_{-1}}\over{\sqrt N}}\right)^2$, 
which clearly has norm $1$. 

\bigskip 
\noindent 
{\bf Lickorish invariant versus Reshetikhin - Turaev invariant} 

In the process of constructing the Lickorish invariant ${\cal F}(M_L)$, 
we have already noticed many similarities between this invariant 
and that of \cite{osp} obtained following a modified Reshetikhin - Turaev 
approach.  Now we prove that these two invariants are actually equivalent, 
namely, on any compact closed orientable $3$ - manifold, 
both invariants take the same value. 

Consider the oriented ribbon graph corresponding to an oriented framed 
link $L$ with $m$ components $L_i$, $i=1, 2, ..., m$.  
We colour the annulus associated with $L_i$ by the $\Uq$ module $W_i$, 
for all $i=1, 2, ..., m$, and denote the resultant coloured ribbon graph 
by $\Gamma(L, \{W_1, ..., W_m\})$. Particularly interesting is the case 
when $W_i = V^{\x l_i}$, where $V$ is the vector module of $\Uq$, 
and $0\le l_i\le N-1$. ($V^{\x 0}={\bf C}$). 
As in the last section, we still denote by $\Gamma(L^{\{l_1, ..., l_m\}})_V$ 
the coloured ribbon graph, which arises from the riented framed link 
$L^{\{l_1, ..., l_m\}}$ obtained by extending each component $L_i$ 
of $L$ to a cable of $l_i$ strands with the same orientation, and 
colouring the annulus associated with each strand by the vector module 
$V$ of $\Uq$. Then a moment's thinking reveals that 
\ban 
F(\Gamma(L^{\{l_1, ..., l_m\}})_V)
&=&F(\Gamma(L, \{V^{\x l_1}, ...,V^{\x l_m}\})). 
\nan   
To examine $\Gamma(L, \{V^{\x l_1}, ...,V^{\x l_m}\})$ more closely, 
we cut open one of its components, say, that associated with $L_1$, 
to obtain another coloured ribbon graph which we denote by  
$\Gamma({\check L}, \{V^{\x l_1}, ...,V^{\x l_m}\})$. Then the  
module homomorphism 
$F(\Gamma({\check L}, \{V^{\x l_1}, ...,V^{\x l_m}\})):$ $V^{\x l_1}$ 
$\rightarrow$ $V^{\x l_1}$ satisfies 
\ban 
F(\Gamma(L, \{V^{\x l_1}, ...,V^{\x l_m}\}))
&=& Str_{V^{\x l_1}} \left[ K F(\Gamma({\check L}, 
\{V^{\x l_1}, ...,V^{\x l_m}\})) \right].\nan   
The right hand side can be evaluated by first decomposing the tensor 
product module $V^{\x l_1}$ into a direct sum of indecomposable 
$\Uq$ modules, then taking the $q$ - supertrace on each module separately. 
Since $l_1\le N-1$, $V^{\x l_1}$ is in fact completely reducible. 
Taking into account the parity of each irreducible submodule ( i.e., the 
evenness or oddness of the highest weight vector), we have 
\ban
& &Str_{V^{\x l_1}} \left[ K F(\Gamma({\check L},
\{V^{\x l_1}, ...,V^{\x l_m}\})) \right]\\  
&=& \sum_{\mu=0}^{N-1} b_\mu^{(l_1)} 
Str_{V^+(\mu) } \left[ K F(\Gamma({\check L},
\{V^+(\mu), V^{\x l_2}, ...,V^{\x l_m}\})) \right]\\ 
&=& \sum_{\mu=0}^{N-1} b_\mu^{(l_1)} F(\Gamma(L,
\{V^+(\mu), V^{\x l_2}, ...,V^{\x l_m}\})),  
\nan
where $ F(\Gamma(L, \{V^+(\mu), V^{\x l_2}, ...,V^{\x l_m}\}))$ 
represents the coloured ribbon graph arising from the oriented 
framed link $L$ with the first component coloured by $V^+(\mu)$, 
and the rest by $V^{\x l_2}$, ..., $V^{\x l_m}$ respectively.  
It can be further expanded by decomposing $V^{\x l_2}$ etc. 
into direct sums of indecomposables.  When $0\le l_i\le N-1$, for all 
$i=1, 2, ..., m$, we obtain the following important relation  
\ba F(\Gamma(L, \{V^{\x l_1}, ...,V^{\x l_m}\}))
=\sum_{\mu_1, ..., \mu_m=0}^{N-1} \  \prod_{i=1}^m b^{(l_i)}_{\mu_i}
\  F(\Gamma(L, \{V^+(\mu_1), ...,V^+(\mu_m)\})). \label{relation}\na

Let us  apply (\ref{relation}) to the definition (\ref{Sigma}) of 
$\Sigma(L)$.  Recalling that $b = d B$, 
we immediately arrive at 
\ban \Sigma(L) &=&\sum_{\mu_1, ..., \mu_m=0}^{N-1}  \  \prod_{i=1}^m b_{\mu_i}
\  F(\Gamma(L, \{V^+(\mu_1), ...,V^+(\mu_m)\})),  \nan 
where the right hand side is precisely what appeared in \cite{osp}. 
Since the $z$ defined by equation (\ref{z}) coincides with the 
corresponding quantity there as well, we easily see that ${\cal F}(M_L)$ 
and the invariant of \cite{osp} are indeed the same.  \\  

\small  
 
\end{document}